\documentclass[11pt]{amsart}
\usepackage{amsmath,amsfonts}
\usepackage[all]{xy}
\newcommand{\Hilb}{\operatorname{Hilb}}
\newcommand{\Pic}{\operatorname{Pic}}
\newcommand{\Hom}{\operatorname{Hom}}

\newcommand{\bul}{\bullet}
\newcommand{\ov}{\mathcal{O}_V}
\newcommand{\oc}{\mathcal{O}_C}
\newcommand{\kv}{\mathcal{K}_V}
\newcommand{\nc}{\mathbb{C}}
\newcommand{\nd}{\mathbb{D}}
\newcommand{\ndtil}{\tilde{\nd}}
\newcommand{\rhotil}{\tilde{\rho}}
\newcommand{\nl}{\mathbb{L}}
\newcommand{\np}{\mathbb{P}}
\newcommand{\nr}{\mathbb{R}}
\newcommand{\nz}{\mathbb{Z}}
\newcommand{\cl}{\mathcal{L}}
\newcommand{\co}{\mathcal{O}}
\newcommand{\pitil}{\tilde{\pi}}
\newcommand{\LG}{\mathcal{L}^\bul}
\newcommand{\spin}{Spin^c(4)}
\newcommand{\spim}{Spin^c(M)}
\newcommand{\spiv}{Spin^c(V)}
	\newcommand{\oo}{{\scriptstyle{\mathcal{O}}}}
	
	\newcommand{\ooo}{{\scriptscriptstyle{\mathcal{O}}}}
	
	\newcommand{\oov}{({\scriptstyle{\mathcal{O}}}_1,{\bf H}_0)}
	
	\newcommand{\ooov}{({\scriptscriptstyle{\mathcal{O}}}_1,{\bf H}_0)}

	\newcommand{\swv}[1]{SW_{V,\ooo}(#1)}
	\newcommand{\swm}[1]{SW_{M,\ooo}(#1)}
	
	\newcommand{\sww}[2]{SW^{#1}_{X,\ooov}(#2)}
	\newcommand{\swwv}[2]{SW^{#1}_{V,\ooov}(#2)}

\newcommand{\lb}{[[}
\newcommand{\rb}{]]}

\newtheorem{thm}{Theorem}
\newtheorem{lem}[thm]{Lemma}
\newtheorem{prop}[thm]{Proposition}
\newtheorem{cor}[thm]{Corollary}
\newtheorem{conj}[thm]{Conjecture}

\theoremstyle{definition}

\newtheorem{dfn}[thm]{Definition}

\newtheorem{rem}[thm]{Remark}

\title[Relations for virtual fundamental classes of Hilbert schemes]{Relations for virtual fundamental classes of Hilbert schemes of curves on surfaces}
\author[M.~D\"urr]{Markus D\"urr$^*$}
\author[Ch.~Okonek]{Christian Okonek$^*$}
\thanks{$^*$Partially supported by: EAGER -- European Algebraic Geometry Research Training Network, contract No. HPRN-CT-2000-00099 (BBW 99.0030), and by SNF, nr. 2000-055290.98/1}
\email{mduerr@math.unizh.ch, okonek@math.unizh.ch}
\begin{document}
\begin{abstract}
In \cite{dko} we constructed virtual fundamental classes $\lb Hilb^m_V \rb$ for Hilbert schemes of divisors of topological type $m$ on a surface $V$, and used these classes to define the Poincar\'e invariant of $V$:
\[
(P^+_V,P^-_V): H^2(V,\nz) \longrightarrow \Lambda^* H^1(V, \nz) \times \Lambda^* H^1(V, \nz)
\]
We conjecture that this invariant coincides with the full Seiberg-Witten invariant computed with respect to the canonical orientation data.

In this note we prove that the existence of an integral curve $C \subset V$ induces relations between some of these virtual fundamental classes $\lb Hilb^m_V \rb$. The corresponding relations for the Poincar\'e invariant can be considered as algebraic analoga of the fundamental relations obtained in \cite{os}. 
\end{abstract}
\maketitle
\section{Introduction}
The symplectic Thom conjecture for homology classes with negative self-intersection, proven by Ozsv\'ath and Szab\'o, is an immediate consequence of the following two facts:
\begin{itemize}
	\item[i)] Taubes' constraints for the Seiberg-Witten basic classes of a closed symplectic four-manifold \cite{t}.
	\item[ii)] A fundamental relation between certain Seiberg-Witten invariants, which arises from embedded surfaces with negative self-inter\-section, due to Ozsv\'ath and Szab\'o \cite{os}.
\end{itemize}
In this note we prove an analoguous relation for the virtual fundamental classes of certain Hilbert schemes of algebraic curves on smooth projective surfaces. To be more precise: Let $V$ be a smooth connected projective surface over $\nc$. For any class $m \in H^2(V,\nz )$ we have the Hilbert scheme $\Hilb^m_V$ parametrizing effective divisors $D \subset V$ with $c_1(\ov (D))=m$. In \cite{dko} we constructed a virtual fundamental class $\lb \Hilb^m_V \rb \in A_*(\Hilb^m_V)$ in the Chow group of $\Hilb^m_V$. Note that there exists a natural morphism $\rho: \Hilb^m_V \to \Pic^m_V$ sending a divisor $D \subset V$ to the class $[\ov(D)]$ of its associated line bundle. Let $\nd  \subset \Hilb^m_V \times V$ be the universal divisor, and put $u:= c_1(\ov (\nd)|_{\Hilb^m_V \times \{ p \}})$, where $p \in V$ is an arbitrary point.

Consider now an integral curve $C \subset V$, set $c:=c_1(\ov(C))$, and denote by $\kappa_c \in \Lambda^2 H^1(V, \nz )^\vee$ the map:
\begin{eqnarray*}
\kappa_c : \Lambda^2 H^1(V, \nz) & \longrightarrow & \nz\\
a \wedge b & \longmapsto & \langle a \cup b \cup c, [V] \rangle.
\end{eqnarray*}
Let $\iota: \Hilb^{m-c}_V \to \Hilb^m_V$ be the closed embedding sending $D' \in \Hilb^{m-c}_V$ to $D'+C \in \Hilb^m_V$. Our main result relates $\lb \Hilb^m_V \rb$ and $\lb \Hilb^{m-c}_V \rb$ when $m \cdot c <0$, and $\lb \Hilb^m_V \rb$ and $\lb \Hilb^{m+c}_V\rb $ when $(k-m)\cdot c <0$. Here $k:=c_1(\kv)$ is the first Chern class of the canonical line bundle.
\setcounter{thm}{2}
\begin{thm}
Let $V$ be a surface, and fix a class $m \in H^2(V,\nz)$. Let $C \subset V$ be a reduced and irreducible curve, and set $c:=c_1(\ov (C))$.
\begin{itemize}
	\item[i)] Suppose that $m\cdot c < 0$, and denote by $\rho$ the map $\Hilb^m_V \to \Pic^m_V$. Let $\iota: \Hilb^{m-c}_V \to \Hilb^m_V$ be the inclusion given by the addition $D \mapsto D + C$. Then we have
\[
\lb \Hilb^m_V \rb = \left( \sum_i \rho^* \left( \frac{\kappa_c^i}{i!}\right) \cdot u^{\frac{c^2 + c \cdot m}{2}-m \cdot c-i}\right) \cap \iota_* \lb \Hilb^{m-c} \rb.
\]
	\item[ii)] Suppose that $(k-m)\cdot c < 0$, and denote by $\rhotil$ the map $\Hilb^{m+c}_V \to \Pic^{m+c}_V$. Let $\iota: \Hilb^{m}_V \to \Hilb^{m+c}_V$ be the inclusion given by the addition $D \mapsto D + C$. Then we have
\[
\iota_* \lb \Hilb^m_V \rb = \left( \sum_i \rhotil^* \left( \frac{(-\kappa_c)^i}{i!}\right) \cdot u^{\frac{c^2+c \cdot k}{2}-(k-m)c-i} \right) \cap \lb \Hilb^{m+c}_V \rb.
\]
\end{itemize}
\end{thm}

In \cite{dko} we used the virtual fundamental classes $\lb \Hilb^m_V\rb$ to define a map
\[
(P^+_V,P^-_V): H^2(V,\nz) \longrightarrow \Lambda^* H^1(V, \nz) \times \Lambda^* H^1(V, \nz)
\]
which we call the Poincar\'e invariant of $V$. This map is invariant under smooth deformations of $V$, satisfies a blow-up formula, and a wall crossing formula for surfaces with $p_g(V)=0$. We conjecture that the Poincar\'e invariant coincides with the full Seiberg-Witten invariant of \cite{ot} computed with respect to the canonical orientation data. Our relations between the virtual fundamental classes of Hilbert schemes lead to corresponding relations for the Poincar\'e invariant:
\setcounter{thm}{5}
\begin{thm}
Let $V$ be a surface, and fix a class $m \in H^2(V,\nz)$. Let $C \subset V$ be a reduced and irreducible curve, and set $c:=c_1(\ov (C))$.
\begin{itemize}
	\item[i)] If $m\cdot c < 0$, then
\[
P^\pm_V(m)=  \tau_{m(m-k)}\left( \exp (\kappa_c)\cap P^\pm_V (m-c) \right).
\]
	\item[ii)] If $(k-m) \cdot c < 0$, then
\[
P^\pm_V(m) = \tau_{m(m-k)} \left( \exp (-\kappa_c) \cap P^\pm_V (m+c) \right).
\]
\end{itemize}
\end{thm}
This result can be considered as an algebraic analog of the Ozsv\'ath-Szab\'o relation, as we will explain in the section 4 below.
\setcounter{thm}{0}
\section{Comparing virtual fundamental classes of Hilbert schemes}
In this paper all surfaces will be smooth, projective, connected, and defined over the field of complex numbers. We denote by $k:=c_1(\kv)$ the first Chern class of the canonical line bundle of a surface $V$.

Recall that an element $c \in H^2(V,\mathbb{Z})$ is {\em characteristic} iff $c \equiv k \mod{2}$. For a characteristic element $c \in H^2(V,\mathbb{Z})$, we denote by $\theta_c \in \Lambda^2 H^1(V,\mathbb{Z})^\vee$ the map
\begin{eqnarray*}
\theta_c: \Lambda^2 H^1(V,\mathbb{Z}) & \longrightarrow & \mathbb{Z}\\
a \wedge b & \longmapsto & \frac{1}{2} \langle a \cup b \cup c, [V]\rangle .
\end{eqnarray*}
We define $\xi_V \in \Lambda^4 H^1(V,\mathbb{Z})^\vee$ to be the map
\begin{eqnarray*}
\xi_V: \Lambda^4 H^1(V,\mathbb{Z}) & \longrightarrow & \mathbb{Z}\\
a \wedge b \wedge c \wedge d & \longmapsto & \langle a \cup b \cup c \cup d, [V]\rangle .
\end{eqnarray*}
\begin{lem}\label{lem:chernchar}
Let $V$ be a surface, and fix a class $m \in NS(V)$. Choose a normalized Poincar\'e line bundle $\nl$ on $\Pic^m_V \times V$, and let $\mu: \Pic^m_V \times V \to \Pic^m_V$ be the projection. Then we have
\[
ch(\mu_! \nl) = \chi (\ov) + \frac{m(m-k)}{2}-\theta_{2m-k}+ \xi_V.
\]
\end{lem}
\begin{proof}
By the Grothendieck-Riemann-Roch theorem \cite[Thm.15.2]{Fu} we have
\[
td(\Pic^m_V) \cdot ch (\mu_! \nl) = \mu_! \left\{ td (\Pic^m_V \times V) \cdot ch (\nl) \right\}.
\]
Hence we need to compute those components of the expression
\[
td (\Pic^m_V \times V) \cdot ch (\nl)
\]
which have bidegree $(*,4)$ with respect to the decomposition
\begin{eqnarray*}
H^*(\Pic^m_V \times V,\mathbb{Z}) & \cong & H^*(\Pic^m_V,\mathbb{Z}) \otimes H^*(V,\mathbb{Z})\\
& \cong & \Lambda^* H^1(V,\mathbb{Z})^\vee \otimes H^*(V,\mathbb{Z}).
\end{eqnarray*}
Set $f:=c_1(\nl)$. Then
\begin{eqnarray*}
f^{2,0} &=& 0 \in H^2(\Pic^m_V,\mathbb{Z}),\\
f^{1,1} &=& id \in \Hom (H^1(V,\mathbb{Z}),H^1(V,\mathbb{Z})),\\
f^{0,2} &=& m \in H^2(V,\mathbb{Z}),
\end{eqnarray*}
where the first equality holds since $\nl$ is normalized.

Next we compute $g:=f^2$. We obtain
\begin{eqnarray*}
g^{2,2} &=& -2\cdot (a \wedge b \mapsto a\cup b ) \in \Hom (\Lambda^2 H^1(V,\mathbb{Z}), H^2(V,\mathbb{Z})),\\
g^{1,3} &=& 2\cdot (a \mapsto a \cup m) \in \Hom (H^1(V,\mathbb{Z}),H^3(V,\mathbb{Z})),\\
g^{0,4} &=& m\cup m \in H^4(V,\mathbb{Z}),
\end{eqnarray*}
all other components being zero. Here the first equality needs justification. Choose a basis $v_1,\ldots ,v_{2q}$ of $H^1(V,\mathbb{Z})$, and denote by $w_1,\ldots , w_{2q}$ the dual basis of $H^1(V,\mathbb{Z})^\vee$. Then
\[
f^{1,1}= \sum_i w_i\otimes v_i,
\]
and
\begin{eqnarray*}
g^{2,2} &=& \left( f^{1,1} \right)^2\\
&=& (\sum_i w_i\otimes v_i)\cup (\sum_i w_i\otimes v_i)\\
&=& - \sum_i\sum_j (w_i\wedge w_j) \otimes (v_i \cup v_j)\\
&=& -2 \sum_{i<j} (w_i\wedge w_j) \otimes (v_i \cup v_j).
\end{eqnarray*}
Now we compute the component of $f^3$ of bidegree $(2,4)$, the only component that does not vanish. We find
\begin{eqnarray*}
f^3 &=& 3 (f^{1,1})^2 \cup f^{0,2}\\
&=& -6 \cdot ( a \wedge b \mapsto a \cup b \cup m) \in \Hom ( \Lambda^2 H^1(V,\mathbb{Z}),H^4(V,\mathbb{Z})).
\end{eqnarray*}
Finally we obtain
\begin{eqnarray*}
f^4 &=& (f^{1,1})^4\\
&=& \sum_{i,j,k,l} (w_i \wedge w_j \wedge w_k \wedge w_l) \otimes ( v_i \cup v_j \cup v_k \cup v_l)\\
&=& 24 \left(\sum_{i< j < k < l} (w_i \wedge w_j \wedge w_k \wedge w_l) \otimes ( v_i \cup v_j \cup v_k \cup v_l) \right)\\
&=& 24 (a \wedge b \wedge c \wedge d \mapsto a \cup b \cup c \cup d).
\end{eqnarray*}
Since $td(\Pic^m_V)=1$, we get
\begin{eqnarray*}
td(\Pic^m_V \times V) &=& pr_V ^* td(V)\\
&=& pr_V^*(1 - \frac{1}{2} k + \chi (\ov ) \cdot PD [pt] ),
\end{eqnarray*}
where $pr_V: \Pic^m_V \times V \to V$ denotes the projection onto $V$.

Putting everything together, we get
\begin{eqnarray*}
ch (\mu_! \nl ) &=& \left\{ \exp f  \cup pr_V^* \left( 1 -\frac{k}{2} + \chi (\ov) \cdot PD[pt] \right) \right\} /[V]\\
&=& \left\{ (\exp f)^{*,4} - (\exp f)^{*,2}\cup pr_V^* \frac{k}{2}+\chi (\ov) \cdot PD[pt] \right\} /[V]\\
&=& \chi (\ov ) + \frac{m \cdot (m-k)}{2}- \theta_{2m-k}+\xi_V.
\end{eqnarray*}
\end{proof}
For an arbitrary element $c \in H^2(V,\nz)$, we denote by $\kappa_c \in \Lambda^2 H^1(V,\nz )^\vee$ the map
\begin{eqnarray*}
\kappa_c: \Lambda^2 H^1(V,\nz) &  \longrightarrow & \nz\\
a \wedge b & \longmapsto & \langle a \cup b \cup c,[V] \rangle.
\end{eqnarray*}
\begin{cor}\label{cor:chern}
Let $V$ be a surface, and fix two classes $m,c \in NS(V)$. Choose a normalized Poincar\'e line bundle $\nl$ on $\Pic^m_V \times V$ and a line bundle $\cl_c$ on $V$ with $c_1(\cl_c)=c$. Let $\mu: \Pic^m_V \times V \to \Pic^m_V$ and $pr_V: \Pic^m_V \times V \to V$ be the projections. Then
\begin{eqnarray*}
ch(\mu_! \nl - \mu_! (\nl \otimes pr_V^* \cl_c^\vee)) &=& m \cdot c - \frac{c^2 + c \cdot k}{2}-\kappa_c,\\
c(\mu_! \nl - \mu_! (\nl \otimes pr_V^* \cl_c^\vee)) &=& \exp (-\kappa_c).
\end{eqnarray*}
\end{cor}
\begin{proof}
The assertion concerning the Chern character is a direct consequence of Lemma \ref{lem:chernchar}. The formula for the Chern class follows immediately since $H^*(\Pic^m_V,\mathbb{Z})$ has no torsion. 
\end{proof}
In order to state our main result, we have to recall some facts from \cite{dko}.

For a surface $V$ and a class $m \in H^2(V,\nz)$, we denote by $\Hilb^m_V$ the Hilbert scheme of divisors $D$ with $c_1(\ov(D))=m$. Let $\nd \subset \Hilb^m_V \times V$ be the universal divisor, and denote by $\pi: \Hilb^m_V \times V \to \Hilb^m_V$ the projection onto $\Hilb^m_V$.

In \cite{dko}, we constructed an obstruction theory (in the sense of Behrend and Fantechi)
\[
\varphi: (R^\bul \pi_* \co_\nd (\nd ))^\vee \to \LG_{\Hilb^m_V}
\]
for $\Hilb^m_V$, and showed that this obstruction theory defines a virtual fundamental class
\[
\lb \Hilb^m_V \rb \in A_{\frac{m(m-k)}{2}} (\Hilb^m_V).
\]  
Choose a point $p \in V$ and set
\[
u:= c_1 (\co(\nd)|_{\Hilb^m_V \times \{ p \}}).
\]
\begin{thm}\label{thm:vfc}
Let $V$ be a surface, and fix a class $m \in H^2(V,\nz)$. Let $C \subset V$ be a reduced and irreducible curve, and set $c:=c_1(\ov (C))$.
\begin{itemize}
	\item[i)] Suppose that $m\cdot c < 0$, and denote by $\rho$ the map $\Hilb^m_V \to \Pic^m_V$. Let $\iota: \Hilb^{m-c}_V \to \Hilb^m_V$ be the inclusion given by the addition $D \mapsto D + C$. Then we have
\[
\lb \Hilb^m_V \rb = \left( \sum_i \rho^* \left( \frac{\kappa_c^i}{i!}\right) \cdot u^{\frac{c^2 + c \cdot m}{2}-m \cdot c-i}\right) \cap \iota_* \lb \Hilb^{m-c} \rb.
\]
	\item[ii)] Suppose that $(k-m)\cdot c < 0$, and denote by $\rhotil$ the map $\Hilb^{m+c}_V \to \Pic^{m+c}_V$. Let $\iota: \Hilb^{m}_V \to \Hilb^{m+c}_V$ be the inclusion given by the addition $D \mapsto D + C$. Then we have
\[
\iota_* \lb \Hilb^m_V \rb = \left( \sum_i \rhotil^* \left( \frac{(-\kappa_c)^i}{i!}\right) \cdot u^{\frac{c^2+c \cdot k}{2}-(k-m)c-i} \right) \cap \lb \Hilb^{m+c}_V \rb.
\]
\end{itemize}
\end{thm}
\begin{proof}
Suppose first that $m\cdot c <0$. Then we have $H^0(\oc(D))=0$ for any divisor $D \in \Hilb^m_V$. It follows that the inclusion $\Hilb^{m-c}_V \to Hilb^m_V$ is an isomorphism. However, the obstruction theories differ: Denote by $\nc$ the product $\Hilb^m_V \times C$. The short exact sequence
\[
0 \to \co_{\nd -\nc} (\nd -\nc) \to \co_\nd(\nd) \to \co_{\nc} (\nd) \to 0
\]
gives rise to a distinguished triangle:
\[
\xymatrix{
R^\bul\pi_* \co_{\nd -\nc} (\nd -\nc) \ar[r] & R^\bul\pi_*\co_\nd(\nd) \ar[d]\\
& R^\bul\pi_* \co_{\nc} (\nd) \ar[ul]^{[1]}}
\]
Here $\pi: \Hilb^m_V \times V \to \Hilb^m_V$ is the projection.
By the excess intersection formula \cite[Prop.1.16]{dko}, we have
\[
\lb \Hilb^m_V \rb = c_{top} (R^1\pi_* \co_{\nc} (\nd)) \cap \iota_* \lb \Hilb^{m-c}_V\rb .
\]
The complex $R^\bul\pi_* \co_{\nc} (\nd))$ is the mapping cone of the morphism
\[
R^\bul\pi_* \co(\nd -\nc) \to R^\bul\pi_* \co(\nd).
\]
Fix a normalized Poincar\'e line bundle $\nl$ on $\Pic^m_V \times V$. Using \cite[Lemma 3.15]{dko}, we see that this choice endows $\Hilb^m_V$ with a relatively ample sheaf $\co_\nl(1)$. Furthermore, there exists an isomorphism
\[
\co(\nd) \stackrel{\cong}{\longrightarrow} (\rho \times id_V)^* \nl \otimes \pi^* \co_\nl(1),
\]
and, since $\nl$ is normalized, we have
\[
u=c_1( \co_\nl(1)).
\]
This implies that $R^\bul\pi_* \co_\nc (\nd)$ is the mapping cone of
\[
\rho^* (R^\bul \mu_* (\nl \otimes pr_V^* \ov(-C))) \otimes \co_\nl(1) \to \rho^* ( R^\bul \mu_* \nl)\otimes \co_\nl(1).
\]
Using Cor.~\ref{cor:chern} we conclude
\[
c_{top} (R^1 \pi_* \co_\nc(\nd)) = \sum_i \rho^* \left( \frac{\kappa_c^i}{i!}\right) \cdot u^{\frac{c^2 + c \cdot m}{2}-m \cdot c-i},
\]
which proves part $i)$.

Suppose now that $(k-m)\cdot c<0$. Then we have $H^1(\oc(D))=0$ for any divisor $D \in \Hilb^{m+c}_V$. Denote by $\ndtil \subset \Hilb^{m+c}_V \times V$ the universal divisor, and let $\pitil: \Hilb^{m+c}_V \times V \to \Hilb^{m+c}_V$ be the projection. It follows that the sheaf $R^1\pitil_* \co_{\Hilb^{m+c}_V\times C}(\ndtil)$ vanishes, and that $\pitil_* \co_{\Hilb^{m+c}_V\times C}(\ndtil)$ is locally free. Moreover, $\iota$ induces an isomorphism
\[
\Hilb^m_V \stackrel{\cong}{\longrightarrow} Z(\lambda),
\]
where $\lambda$ is the canonical section in $\pitil_* \co_{\Hilb^{m+c}_V\times C}(\ndtil)$.

The short exact sequence
\[
0 \to \co_\nd(\nd) \to \co_{\nd+\nc}(\nd+\nc) \to \co_\nc (\nd + \nc) \to 0
\]
gives rise to the following distinguished triangle:
\[
\xymatrix{
R^\bul\pi_* \co_{\nd} (\nd) \ar[r] & R^\bul\pi_*\co_{\nd+\nc}(\nd+\nc) \ar[d]\\
& R^\bul\pi_* \co_{\nc} (\nd+\nc) \ar[ul]^{[1]}}
\]
Hence functoriality \cite[Thm.1]{kkp} yields
\[
\iota_* \lb Hilb^m_V \rb = c_{top}( \pitil_* \co_{\Hilb^{m+c}_V \times C}(\ndtil)) \cap \lb \Hilb^{m+c}_V \rb .
\]
Fix again a normalized Poincar\'e line bundle $\nl$ on $\Pic^m_V$. By arguments similar to those of the first part, we see that $R^\bul \pitil_* \co_{\Hilb^{m+c}_V \times C}(\ndtil)$ is the mapping cone of 
\[
\rhotil^*( R^\bul \mu_* \nl ) \otimes \co_{\nl \otimes pr_V^* \ov(C)}(1) \to  \rhotil^* ( R^\bul\mu_* (\nl \otimes pr_V^* \ov(C))\otimes \co_{\nl \otimes pr_V^* \ov(C)}(1).
\]
Now Cor.~\ref{cor:chern} implies
\[
c_{top}( \pitil_* \co_{\Hilb^{m+c}_V \times C}(\ndtil)) =\sum_i \rhotil^* \left( \frac{(-\kappa_c)^i}{i!}\right) \cdot u^{\frac{c^2+c \cdot k}{2}-(k-m)c-i}
\]
\end{proof}
\begin{rem}
When $C$ is rational, i.e.~when the normalization $\hat{C}$ is isomorphic to $\np^1$, then $\kappa_c=0$. When $C$ is isomorphic to $\np^1$ and $c^2 \in \{ 0,-1\}$, then $m \cdot c <0$ or $(k-m) \cdot c <0$ for any $m \in H^2(V,\nz)$. 
\end{rem}
To see this, let $j: \hat{C} \to V$ be the map induced by the inclusion $C \subset V$. Then for all $a,b\in H^1(V,\nz)$
\begin{eqnarray*}
\kappa_c(a\wedge b) &=& \langle a \cup b, j_* [ \hat{C} ] \rangle\\
&=& \langle j^*a \cup j^*b ,[\hat{C}] \rangle.
\end{eqnarray*}
Since the curve $\hat C$ is simply connected, the pull-backs $j^*a$ and $j^*b$ vanish, and therefore
\[
\kappa_c(a\wedge b)=0.
\]
When $C$ is isomorphic to $\np^1$ and $c^2 \in \{ 0,-1\}$, the adjunction formula yields $k \cdot c <0$. This proves the second claim.
\section{Relations for Poincar\'e invariants and the adjunction inequality}
First we recall the definition of the Poincar\'e invariant.
Let $V$ be a surface, $p\in V$ an arbitrary point. Fix a class $m \in H^2(V,\mathbb{Z})$, denote by $\nd^+$ the universal divisor over the Hilbert scheme $\Hilb^m_V$, and set
\[
u^+:= c_1\left( \co (\nd^+)|_{ \Hilb^m_V \times \{ p \}}\right) \in H^2(\Hilb^m_V,\mathbb{Z}).
\]
Since $V$ is connected, the class $u^+$ does not depend on the chosen point $p$. Likewise, denote by $\nd^-$ the universal divisor over the Hilbert scheme $\Hilb^{k-m}_V$, where $k=c_1(\kv)$. Put
\[
u^-:= c_1\left( \co (\nd^-)|_{ \Hilb^{k-m}_V \times \{ p \}}\right) \in H^2(\Hilb^{k-m}_V,\mathbb{Z}).
\]
Denote by $\rho^\pm$ the following morphisms:
\begin{eqnarray*}
\rho^+: \Hilb^m_V & \longrightarrow & \Pic^m_V\\
D & \longmapsto & [\ov (D)]
\end{eqnarray*}
\begin{eqnarray*}
\rho^-: \Hilb^{k-m}_V & \longrightarrow & \Pic^m_V\\
D' & \longmapsto & [\kv (-D')]
\end{eqnarray*}
By abuse of notation, we will denote the image of $\lb \Hilb^m_V \rb$ under the cycle map $A_*(\Hilb^m_V) \to H_*(\Hilb^m_V,\nz )$ by the same symbol.
\begin{dfn}
Let $V$ be a surface.  The \emph{Poincar\'e invariant} of $V$ is the map
\begin{eqnarray*}
(P^+_V,P^-_V) : H^2(V, \mathbb{Z}) & \longrightarrow & \Lambda^* H^1(V,\mathbb{Z}) \times
\Lambda^* H^1(V,\mathbb{Z})\\
m & \longmapsto & (P^+_V(m),P^-_V(m)),
\end{eqnarray*}
defined by
\[
P^+_V(m):= \rho^+_* \left( \sum_i (u^+)^{i} \cap [[Hilb^m_V]] \right)
\]
and
\[
P^-_V(m):= (-1)^{ \chi (\ov ) + \frac{m(m-k)}{2}} \rho^-_* \left( \sum_i
(-u^-)^{i} \cap [[Hilb^{k-m}_V]] \right),
\]
if $m \in NS(V)$, and by $P^\pm_V(m):=0$ otherwise.
\end{dfn}
For an integer $n$ we define a truncation map
\[
\tau_{\leq n} : \Lambda^* H^1(V,\nz ) \longrightarrow \Lambda^* H^1(V,\nz )
\]
as follows: when $P=\sum_i P_i$ is the decomposition of a form $P$ into its homogeneous components $P_i \in \Lambda^i H^1(V,\nz )$, then
\[
\tau_{\leq n}(P):= \sum_{i=0}^n P_i.
\]
\begin{thm}\label{thm:rel}
Let $V$ be a surface, and fix a class $m \in H^2(V,\nz)$. Let $C \subset V$ be a reduced and irreducible curve, and set $c:=c_1(\ov (C))$.
\begin{itemize}
	\item[i)] If $m\cdot c < 0$, then
\[
P^\pm_V(m)=  \tau_{m(m-k)}\left( \exp (\kappa_c)\cap P^\pm_V (m-c) \right).
\]
	\item[ii)] If $(k-m) \cdot c < 0$, then
\[
P^\pm_V(m) = \tau_{m(m-k)} \left( \exp (-\kappa_c) \cap P^\pm_V (m+c) \right).
\]
\end{itemize}
\end{thm}
\begin{proof}
Suppose that $m\cdot c <0$, and let $\iota^+$ be the inclusion $\Hilb^{m-c}_V \to \Hilb^m_V$. By part $i)$ of Thm.~\ref{thm:vfc} we have
\begin{eqnarray*}
P^+_V(m) &=& \rho^+_* \left(\sum_i u^i \cap \lb \Hilb^m_V \rb \right)\\
&=& \rho^+_* \left( \sum_i u^i \cap \left( \sum_j (\rho^+)^* \left( \frac{ \kappa_c^j}{j!} \right) u^{\frac{c^2+c\cdot k}{2} -m\cdot c -j} \right) \cap \iota^+_* \lb \Hilb^{m-c}_V \rb \right)\\
&=&
\sum_j \frac{\kappa_c^j}{j!} \cap \rho^+_* \left( \sum_i u^{i+\frac{c^2+c \cdot m}{2} -m \cdot c -j} \cap \iota^+_* \lb \Hilb^{m-c}_V \rb \right)\\
&=& \tau_{m(m-k)} \left( \exp (\kappa_c) \cap P^+_V(m-c) \right).
\end{eqnarray*}
Let $\iota^-$ be the inclusion $\Hilb^{k-m}_V \to \Hilb^{k-m+c}_V$, and set $\epsilon:= (-1)^{\chi (\ov)+\frac{m(m-k)}{2}}$. Note that under the isomorphism
\begin{eqnarray*}
\Pic^m_V & \longrightarrow & \Pic^{k-m}_V\\
{[\cl ]} & \longmapsto & [\kv \otimes \cl^\vee]
\end{eqnarray*}
the cohomology class $\kappa_c$ is mapped to $\kappa_c$, since this class is of degree 2. Hence part $ii)$ of Thm.~\ref{thm:vfc} yields
\begin{eqnarray*}
P^-_V(m) &=& \epsilon \cdot (\rho^-)_* \left( \sum_i (-u)^i \cap \iota^-_* \lb \Hilb^{k-m}_V \rb \right)\\
&=& \epsilon \cdot \rho^-_* \left( \sum_i (-u)^i \cap \left( \sum_j (\rho^-)^* \left( \frac{ (-\kappa_c)^j}{j!} \right) \cdot u^{\frac{c^2+c\cdot k}{2} -m\cdot c -j} \cap \lb \Hilb^{k-m+c}_V \rb \right) \right)\\
&=& \epsilon \cdot (-1)^{\frac{c^2+c\cdot k}{2} -m\cdot c}\left( \sum_j \frac{\kappa_c^j}{j!} \cap \rho^-_* \left( \sum_i (-u)^{i+\frac{c^2+c\cdot k}{2} -m\cdot c -j} \cap \lb \Hilb^{k-m+c}_V \rb \right) \right)\\
&=& \tau_{m(m-k)} (\exp (\kappa_c) \cap P^-_V(m-c)).
\end{eqnarray*}
The proof in the case $(k-m)\cdot c <0$ is similar. We omit the details.
\end{proof}
Recall that a class $m \in H^2(V,\nz)$ is basic for a surface $V$, if
\[
(P^+_V(m),P^-_V(m)) \neq (0,0).
\]
The surface $V$ is of simple type if all basic classes $m \in H^2(V,\nz )$ satisfy $m(m-k)=0$. In \cite[Prop.6.25]{dko} we have shown that surfaces with $p_g(V)>0$ are of simple type. The following result can be considered as an algebraic analog of the Ozsv\'ath-Szab\'o inequality \cite[Cor.1.7]{os}.
\begin{prop}
Let $V$ be a surface with $p_g(V)>0$, let $C \subset V$ be a curve, and set $c:= c_1(\ov (C))$. For any basic class $m \in H^2(V,\nz)$ we have
\[
0 \leq m \cdot c \leq k \cdot c,
\]
unless $C$ is a smooth rational curve. In this case we have
\[
-1 \leq m \cdot c \leq k \cdot c+1
\]
for all basic classes $m \in H^2(V,\nz)$.
\end{prop}
\begin{proof}
Assume first that $m$ is a basic class with $m\cdot c <0$. Then Thm.~\ref{thm:rel} implies that also $m-c$ is a basic class. We have
\[
\frac{(m-c)(m-c-k)}{2} = \frac{m(m-k)}{2} + p_a(C)-1-m \cdot c
\]
Since any surface $V$ with $p_g(V)>0$ is of simple type, this implies
\[
p_a(C)=0 \text{ and } m\cdot c=-1.
\]
Analoguously, if $m$ is a basic class with $m \cdot c > k \cdot c$, then also $m+c$ is a basic class. Because
\[
\frac{(m+c)(m+c-k)}{2} = \frac{m(m-k)}{2} + p_a(C)-1 - (k-m) \cdot c,
\]
we obtain this time
\[
p_a(C)=0 \text{ and } (k-m)\cdot c=-1.
\]
\end{proof}
\section{Connection with the Ozsv\'ath-Szab\'o relation}
In order to explain the connection between Thm.~\ref{thm:rel} and the Ozsv\'ath-Szab\'o relation, we briefly recall the structure of the full Seiberg-Witten invariants; for the construction and details, we refer to \cite{ot}.

Let $(M,g)$ be a closed oriented Riemannian 4-manifold with first Betti number $b_1$. We denote by $b_+$ the dimension of a maximal subspace of $H^2(M,\nr)$ on which the intersection form is positive definite. Recall that the set of isomorphism classes of $\spin$-structures on $(M,g)$ has the structure of a $H^2(M,\mathbb{Z})$-torsor. This torsor does, up to a canonical isomorphism, not depend on the choice of the metric $g$ and will be denoted by $\spim$.

We have the Chern class mapping
\begin{eqnarray*}
c_1:\spim & \longrightarrow & H^2(M,\mathbb{Z}) \\
\mathfrak{c} & \longmapsto & c_1(\mathfrak{c}),
\end{eqnarray*}
whose image consists of all characteristic elements.

If $b_+>1$, then the Seiberg-Witten invariants are maps
\[
SW_{M,\ooo} : \spim \longrightarrow \Lambda^*H^1(M,\mathbb{Z}),
\]
where $\oo$ is an orientation parameter.

When $b_+=1$, then the invariants depend on a chamber structure and are maps
\[
(SW^+_{M,\ooov},SW^-_{M,\ooov}) : \spim \longrightarrow \Lambda^*H^1(M,\mathbb{Z}) \times\Lambda^*H^1(M,\mathbb{Z}),
\]
where $\oov$ are again orientation data. The difference of the two components is a purely topological invariant.

Let $\Sigma \subset M$ be a smoothly embedded, oriented, closed two-manifold. Fix a standard symplectic basis for $H_1(\Sigma ,\nz)$ and let $\{ A_i,B_i \}_{i=1}^g$ be its image in $H^1(M,\nz)^\vee$. We define the class $\theta(\Sigma) \in \Lambda^2 H^1(M,\nz)^\vee$ by
\[
\theta(\Sigma)= \sum_i A_i \wedge B_i.
\]
\begin{thm}[Ozsv\'ath-Szab\'o]\label{thm:os}
Let $M$ be a closed, oriented, smooth four-manifold with $b_+>0$, and let $\Sigma \subset M$ be a smoothly embedded, oriented, closed two-manifold of genus $g>0$ with negative self-intersection
\[
[\Sigma] \cdot [\Sigma] =-n.
\]
If $b_+>1$, then for each $\spin$-structure $\mathfrak{c}$ with expected dimension $d(\mathfrak{c})\geq 0$ and
\[
| \langle c_1(\mathfrak{c}),[\Sigma] \rangle | \geq 2g+n
\]
we have
\[
\swm{\mathfrak{c}} = \tau_{\leq d(\mathfrak{c})} (\exp (\theta(\epsilon \Sigma))\cap \swm{\mathfrak{c}+ \epsilon PD(\Sigma)}),
\]
where $\epsilon=\pm 1$ is the sign of $\langle c_1(\mathfrak{c}),[\Sigma] \rangle$, and $PD(\Sigma)$ denotes the class Poincar\'e dual to $[\Sigma]$.

If $b_+=1$, then for each $\spin$-structure $\mathfrak{c}$ with expected dimension $d(\mathfrak{c})\geq 0$ and
\[
| \langle c_1(\mathfrak{c}),[\Sigma] \rangle | \geq 2g+n
\]
we have

\[
\sww{\pm}{\mathfrak{c}} = \tau_{\leq d(\mathfrak{c})} (\exp(\theta(\epsilon \Sigma))\cap \sww{\pm}{\mathfrak{c}+ \epsilon PD[\Sigma]}).
\]
\end{thm}
We need the following
\begin{lem}
Let $M$ be a closed, oriented, smooth four-manifold. Let $\Sigma \subset M$ be a smoothly embedded, oriented, closed two-manifold, and let $c$ be the Poincar\'e dual of the homology class $[\Sigma]$. Then
\[
\theta(\Sigma)(a \wedge b) = \langle a \cup b \cup c,[M] \rangle \ \ \forall a,b \in H^1(M,\nz ).
\]
\end{lem}
\begin{proof}
Fix a standard symplectic basis $\{ \alpha_i , \beta_i \}_{i=1}^g$, and let $\{A_i, B_i \}_{i=1}^g$ be its image in $H^1(M,\nz)^\vee$. Then for all $a,b \in H^1(M,\nz)$
\begin{eqnarray*}
\langle a \cup b \cup c,[M] \rangle &=& \langle a \cup b , c \cap [M] \rangle\\
&=& \langle a \cup b , j_* [\Sigma] \rangle\\
&=& \langle j^*a \cup j^*b, [\Sigma ]\rangle\\
&=& \sum_{i=1}^g \det 
\begin{pmatrix} j^*a (\alpha_i) & j^* a (\beta_i)\\
j^*b (\alpha_i) & j^* b (\beta_i)\end{pmatrix}\\
&=& \sum_{i=1}^g \det 
\begin{pmatrix} A_i (a) & B_i(a)\\
A_i(b) & B_i (b)\end{pmatrix}\\
&=& \theta(\Sigma)(a \wedge b).
\end{eqnarray*}
\end{proof}
At this point it is clear, that Thm.~\ref{thm:rel} and Thm.~\ref{thm:os} are fomally analoguous statements. We believe however, that the actual source of this analogy is the conjectured equivalence between our Poincar\'e invariants and the full Seiberg-Witten invariants.
To be precise, let $V$ be a surface. Any {\em Hermitian} metric $g$ on $V$ defines a {\em canonical} $\spin$-structure on $(V,g)$. Its class $\mathfrak{c}_{can} \in \spiv$ does not depend on the choice of the metric. The Chern class of $\mathfrak{c}_{can}$ is $c_1(\mathfrak{c}_{can})=-c_1(\kv )=-k$.

Since $\spiv$ is a $H^2(V,\nz)$-torsor, the distinguished element $\mathfrak{c}_{can}$ defines a bijection:
\begin{eqnarray*}
H^2(V,\nz) & \longrightarrow & \spiv\\
m & \longmapsto & \mathfrak{c}_m
\end{eqnarray*}
The Chern class of the twisted structure $\mathfrak{c}_m$ is $2m-k$. Recall that any surface defines canonical orientation data $\oo$ and $\oov$ respectively.

The precise conjectured relation between Poincar\'e and Seiberg-Witten invariants is:
\begin{conj}\label{conj:main}
Let $V$ be a surface, and denote by $\oo$ or $\oov$ the canonical orientation data. If $p_g(V)=0$, then
\[
P^\pm_V(m)=\swwv{\pm}{\mathfrak{c}_m}\ \ \forall m\in H^2(V,\mathbb{Z}).
\]
If $p_g(V)>0$, then
\[
P^+_V(m)=P^-_V(m)=\swv{\mathfrak{c}_m}\ \ \forall m\in H^2(V,\mathbb{Z}).
\]
\end{conj}
If this conjecture holds, Thm.~\ref{thm:rel} is essentially a consequence of Thm.~\ref{thm:os}. 
To see this, let $C \subset V$ be an integral curve in the surface $V$. Its arithmetic genus is given by the adjunction formula
\[
p_a(C)= \frac{c^2+c\cdot k}{2} +1,
\]
where $c:=c_1(\ov(C))$. Hence the inequality
\[
| \langle c_1(\mathfrak{c}),[\Sigma] \rangle | \geq 2g+n
\]
with $n=-[\Sigma]\cdot [\Sigma]$ reads
\[
| \langle c_1(\mathfrak{c}),[\Sigma] \rangle | \geq c \cdot k+2.
\]
When $\mathfrak{c}=\mathfrak{c}_m$ for some $m \in H^2(V,\nz )$, this means
\[
| (2m-k)\cdot c | \geq c\cdot k +2,
\]
or equivalently
\[
m \cdot c \leq -1 \text{ or } (k-m)\cdot c \leq -1.
\]
Moreover, in the first case $\epsilon =-1$, whereas in the second case $\epsilon=+1$.

Conversely, Thm.~\ref{thm:rel} yields further evidence for the truth of Conj.~\ref{conj:main}.

\end{document}